\documentclass[letterpaper, 10 pt, conference]{ieeeconf}  

\IEEEoverridecommandlockouts                              

\usepackage{tikz}
\usetikzlibrary{calc} 
\usepackage{macros}
         \usepackage{epsfig} 
         \usepackage{graphicx}
         \usepackage{graphics}
         \usepackage{latexsym} 
         \usepackage{amsmath}
         \usepackage{psfrag}
         \usepackage{amssymb}
         \usepackage{multicol}
         \usepackage{cite}
         \usepackage{float}
         
\title{\LARGE \bf
Circular orbit spacecraft control at the L4 point using Lyapunov functions
}

\author{R Agrawal$^{1}$ and R N Banavar$^{2}$
\thanks{$^{1}$Rachana Agrawal is with the Department of Aerospace Engineering, Indian Institute of Technology Bombay}
\thanks{$^{2}$Ravi N Banavar is with Systems and Control
Engineering, Indian Institue of Technology Bombay}
}

\begin{document}

\maketitle
\thispagestyle{empty}
\pagestyle{empty}

\begin{abstract}

The objective of this work is to demonstrate the utility of Lyapunov functions in control synthesis
for the purpose of maintaining and stabilizing a spacecraft in a circular orbit at the L4 point in the circular restricted three body problem (CR3BP). Incorporating the requirements of a fixed radius orbit and a desired angular momentum, a Lyapunov function is constructed and the requisite analysis is performed to obtain a controller. Asymptotic stability is proved in a defined region around the L4 point with the help of LaSalle's principle. 

\end{abstract}

\section{INTRODUCTION}
The restricted three body problem is defined as follows: Two bodies revolve around their
center of mass in circular orbits under the influence of their mutual gravitational attraction
and a third body (attracted by the previous two but not influencing their motion) moves in the
plane defined by the two revolving bodies. The restricted problem of three bodies is to describe
the motion of this third body \cite{szebehely1969theory}.
\par 
Libration points are the natural equilibrium solutions of the restricted three-body problem (R3BP).In the last few years, the interest concerning the libration points for space applications has
risen within the scientific community ~\cite{canalias2004assessment}. This is because the libration points offer the unique possibility to have a fixed configuration with respect to two primaries. 
Therefore, a libration point mission could fulfill a lot of mission constraints that are not achievable with the classical Keplerian two-body orbits. Moreover, exploiting the stable and unstable part of the dynamics
related to these equilibria, low-energy station keeping missions of practical interest can be
obtained. 
\par 
Considerable work has been done on the CR3BP with linearized dynamics \cite{ADAB1967stability} \cite{AM1962onthe} \cite{markeev1969stability}. Work on existence of formation flight trajectories near the triangular libration point in the CR3BP using the linearized equations of motion has been presented in \cite{catlin2007earth}. Feedback linearization techniques have been applied to various formation flights near the vicinity of the libration points \cite{vadali2004design} \cite{marchand2005control}. An extensive bibliographical survey of problems in this context can be found in \cite{gómez2001dynamics}. 
\par 
In \cite{chang2002lyapunov}, a Lyapunov based control is derived to achieve transfer between elliptic Keplerian orbits. Taking inspiration from \cite{chang2002lyapunov}, we attempt to exploit the theory of Lyapunov functions and related stability notions to derive a feedback control law to keep the spacecraft in a circular orbit around the libration point. We use Newton's universal law of gravitation to express the spacecraft dynamics. This approach of using the original nonlinear dynamics of the system to derive the feedback controller has not been adopted in the area of control related to CR3BP.
\par 
The document is organized as follows. In the second section we introduce the notations used in the text and derive the dynamics of the spacecraft in the CR3BP. In the third section we define the characteristics of the desired orbit and propose a candidate Lyapunov function. We carry out further computation and analysis to 
construct a suitable controller so as to render the desired orbit stable. In the fourth section we use LaSalle's invariance principle to prove asymptotic stability of the desired orbit. Finally, in the fifth section, we present the numerical simulation results with respect to the Earth-Moon system. 

\section{Notation and system modelling}
	We first present the notation used in the document and follow this up with a dynamic
	model. 
	\begin{itemize}
\item $m_1, m_2, m_s$ - two primaries and spacecraft
\item $m_1, m_2, L4$ form a plane of rotation - synodic frame
\item $r_{mj}$ - distance of $m_j$ from centre of mass (COM) of $m_1 m_2$
\item $\boldsymbol{r_c}$ - position vector of L4 point from the COM
\item $\boldsymbol{r_{cs}}$ - position vector of spacecraft from L4
\item $\boldsymbol{r_{js}}$ - position vector of spacecraft from $m_j$
\item $\phi$ - angle of rotation of the synodic frame with respect to the sidereal frame
\item subscript $b$ - vectors in synodic frame
\item subscript $i$ - vectors in inertial frame
\item $x_i, y_i$ are along the inertial frame
\item $x_s, y_s$ are along the synodic frame
\item $\boldsymbol{\omega}$ - angular velocity vector of the two primaries about the COM
\item $\hat{\omega}$ - skew symmetric matrix of $\omega$ given by 
\[  \hat{\omega} =  \dot{\phi} \pmat{0 & 1 & 0  \\ -1 & 0  & 0  \\  0 & 0 & 0  }  
	  \]
\end{itemize}

\begin{figure}[h]
\begin{center}
\begin{tikzpicture}[>=latex]
\coordinate (O) at (0,0);
\coordinate (m_s) at ([shift={(120:3cm)}]O);
\coordinate (m_1) at ([shift={(15:-3cm)}]O);
\coordinate (m_2) at ([shift={(15:3cm)}]O);
\coordinate (L4) at ([shift={(105:2.6cm)}]O);
\draw[black,fill=black] ([shift={(15:-3cm)}]O) circle (1.0ex);
\draw[black,fill=black] ([shift={(15:3cm)}]O) circle (1.0ex);
\draw[black,fill=black] ([shift={(120:3cm)}]O) circle (0.5ex);
\foreach \Valor/\p in {m_s/s,m_1/m_1,m_2/m_2,L4/c}
  \draw[->] (O) -- (\Valor) node[midway,above,sloped] {$r_{\p}$};
  \draw[->] (m_1) -- (m_s)
  node[midway,above,sloped] {$r_{1s}$};
  \draw[->] (m_2) -- (m_s)
  node[midway,above,sloped] {$r_{2s}$};
  \draw[->] (L4) -- (m_s)
  node[midway,above,sloped] {$r_{cs}$};
\foreach \Valor/\Pos in {O/below,m_2/right,m_1/left,m_s/left,L4/right}
  \node[\Pos] at (\Valor) {$\Valor$};    
\draw[gray!90] 
  (O) -- ++(0pt,5cm) node[above] {$y_i$};  
\draw[gray!90] 
  ([xshift=-3cm]O) -- ([xshift=3cm]O) node[right] {$x_i$}; 
\draw[dashed]
  (O) -- ([shift={(105:5cm)}]O) node[above] {$y_b$};
\draw[dashed]
  ([shift={(15:-4cm)}]O) -- ([shift={(15:4cm)}]O) node[right] {$x_b$}; 
\end{tikzpicture}
\caption{The three body system with spacecraft in vicinity of L4 point}	
\end{center}
\end{figure}
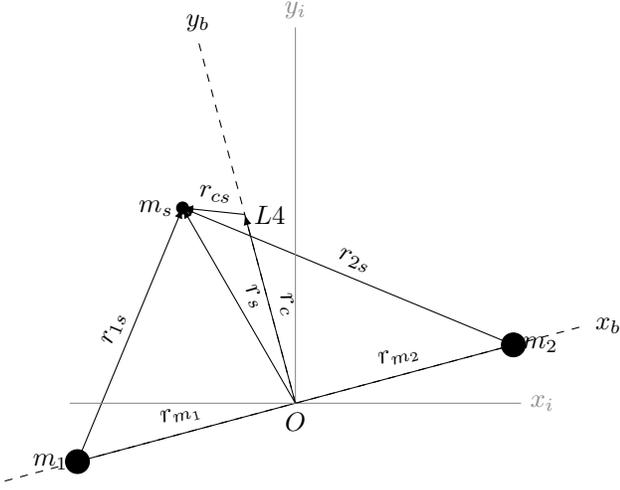
	From Newton's universal law of gravitation, we have the dynamics of the satellite expressed as 
	\begin{equation}
		m_3\boldsymbol{\ddot{r}_{s_i}} = \frac{km_1m_3}{||\boldsymbol{r_{1s_i}}||^2}(-\boldsymbol{\hat{r}_{1s_i}}) + 
			\frac{km_2m_3}{||\boldsymbol{r_{2s_i}}||^2}(-\boldsymbol{\hat{r}_{2s_i}}) 
			\label{newtonlaw}
	\end{equation}
	where $k$ is the gravitational constant, $\boldsymbol{r_{1s_i}}$ and $\boldsymbol{r_{2s_i}}$ are the position vectors of the spacecraft in the inertial frame from masses $m_1$ and $m_2$ respectively. The sidereal (inertial) and
	synodic coordinate representations of the position vectors are related by
	
	\begin{equation}
		\boldsymbol{r_i} = R\boldsymbol{r_b}	\label{synspat}
	\end{equation}
	where $R$ is the rotation matrix relating the synodic frame
		to the sidereal frame, given by the matrix 
	\[ 	\pmat{
		cos\phi & sin\phi & 0 \\
		-sin\phi & cos\phi & 0 \\
		0 & 0 & 1 
		}
	\]
	where $\phi$ denotes the angle of rotation. Note that we are concerned with planar rotations.
	Differentiating (\ref{synspat}) twice, we have 
	\begin{equation}
	\boldsymbol{\ddot{r}_{s_i}} = \ddot{R}\boldsymbol{r_{s_b}} + 2\dot{R}\boldsymbol{\dot{r}_{s_b}} 
	+ R\boldsymbol{\ddot{r}_{s_b}} \label{rddot_i_to_b}
	\end{equation}
	Our objective is to express (\ref{newtonlaw}) in the body frame coordinates. We proceed
	as follows.
We substitute (\ref{rddot_i_to_b}) in (\ref{newtonlaw}). Also,
	\begin{equation*}
	\boldsymbol{r_{1s_i}} = R\boldsymbol{r_{1s_b}} = R(\boldsymbol{r_{1c_b}} + \boldsymbol{r_{cs_b}})
	\end{equation*}
	which yields
	\begin{equation*}
	\begin{split}
		&\frac{km_1m_3}{||\boldsymbol{r_{1s_i}}||^3}(-\boldsymbol{r_{1s_i}}) + 
		\frac{km_2m_3}{||\boldsymbol{r_{2s_i}}||^3}(-\boldsymbol{r_{2s_i}}) \\
			 %
			= &\frac{-Rkm_1m_3}{||\boldsymbol{r_{1s_b}}||^3}(\boldsymbol{r_{1c_b}} + \boldsymbol{r_{cs_b}}) + 
			\frac{-Rkm_2m_3}{||\boldsymbol{r_{2s_i}}||^3}(\boldsymbol{r_{2c_b}} + \boldsymbol{r_{cs_b}})
	\end{split}
	\end{equation*}
	(\ref{newtonlaw}) now reads as
	\begin{equation*}
	\begin{split}
		m_3(\ddot{R}\boldsymbol{r_{s_b}} + 2\dot{R}\boldsymbol{\dot{r}_{s_b}} + R\boldsymbol{\ddot{r}_{s_b}}) 
		= & \frac{-Rkm_1m_3}{||\boldsymbol{r_{1s_b}}||^3}(\boldsymbol{r_{1c_b}} + \boldsymbol{r_{cs_b}}) \\ & + 
			\frac{-Rkm_2m_3}{||r_{2s_S}||^3}(\boldsymbol{r_{2c_b}} + \boldsymbol{r_{cs_b}})
	\end{split}
	\end{equation*}
	
	\noindent
	Regrouping terms%
	\begin{equation*}
	\begin{split}
		\boldsymbol{\ddot{r}_{s_b}} = & -R^{-1}\ddot{R}\boldsymbol{r_{s_b}} - 2R^{-1}\dot{R}
		\boldsymbol{\dot{r}_{s_b}} \\
		& -k(\frac{m_1(\boldsymbol{r_{1c_b}} + \boldsymbol{r_{cs_b}})}{||\boldsymbol{r_{1s_b}}||^3} +
		 \frac{m_2(\boldsymbol{r_{2c_b}} + \boldsymbol{r_{cs_b}})}{||\boldsymbol{r_{2s_b}}||^3})
    \end{split}
	\end{equation*}
	and using the expression
	\begin{equation*} 
		\dot{R} = R\hat{\boldsymbol{\omega}} , \;\;\;\;
		\ddot{R} = \dot{R}\hat{\boldsymbol{\omega}} = R\hat{\boldsymbol{\omega}}\hat{\boldsymbol{\omega}}
	\end{equation*}
	 we have
	 \begin{equation}
	 \begin{split}
		\boldsymbol{\ddot{r}_{cs_b}} = &  \dot{\phi}^2\boldsymbol{r_{cs_b}} + \dot{\phi}^2\boldsymbol{r_{c_b}} -
		 2\hat{\boldsymbol{\omega}}\boldsymbol{\dot{r}_{cs_b}} \\ &- \frac{km_1}{||\boldsymbol{r_{1s_b}}||^3}(\boldsymbol{r_{cs_b}} +
		  \boldsymbol{r_{1c_b}}) - \frac{km_2}{||\boldsymbol{r_{2s_b}}||^3}(\boldsymbol{r_{cs_b}} + \boldsymbol{r_{2c_b}})
		  \end{split} \label{naturaldyn}
	\end{equation}
	\noindent
	Now (\ref{naturaldyn}) represents the natural dynamics of the spacecraft near the L4 point. 
	 Note that for the circular restricted three body problem $\dot{\phi}$ is assumed to be a constant
	  and is given by 
	\begin{equation*}
		\dot{\phi} = \sqrt{\frac{F_{12}}{m_1r_{m_1}}}
		\Rta  \ddot{R} = -\dot{\phi}^2R
	\end{equation*}
	where $F_{12}$ is the gravitational force exerted by $m_2$ on $m_1$.
	Also
\begin{equation*}
\boldsymbol{r_{s_b}} = \boldsymbol{r_{c_b}} + \boldsymbol{r_{cs_b}}  \;\;\;\; \hbox{with} \;\;\;\; \boldsymbol{\ddot{r}_{s_b}} = \boldsymbol{\ddot{r}_{cs_b}}
\end{equation*} 
\noindent
since $\ddot{r}_{c_b} = 0$ as it is the L4 point and is stationary in the synodic frame.
\par 
Now we introduce the control input as follows:
	 \begin{equation}
\boldsymbol{\ddot{r}_{cs_b}} = f(\boldsymbol{r_{cs_b}}, \boldsymbol{\dot{r}_{cs_b}}) + \boldsymbol{\bar{u}_b} \label{acc}
\end{equation}
where $f(\boldsymbol{r_{cs_b}}, \boldsymbol{\dot{r}_{cs_b}})$ is the whole right hand side of (\ref{naturaldyn}),  $\boldsymbol{\bar{u}_b} ( \in \mathbb{R}^3) $ 
denotes the three independent control inputs. Note that in practice control is achieved through force generated by firing thrusters. The control law here has the units of acceleration. The force applied will be acceleration times the mass of spacecraft.  
\section{Lyapunov function and analysis}

	For further analysis, we characterize the orbits around the L4 point by a two-tuple
	$(d, \boldsymbol{L_d})$, where $d$ stands for the desired magnitude of the radius of the orbit and $\boldsymbol{L_d}$
	is the desired angular momentum.
	For a circular orbit of a pre-specified radius around the triangular libration point 
	the candidate Lyapunov function must hence incorporate three objectives:
		\begin{itemize}
		
			\item The velocity vector is perpendicular to the position vector,
			\[ \boldsymbol{r_{cs_b}}\cdot\boldsymbol{\dot{r}_{cs_b}} =   0  \]
			\item The angular momentum is constant,
			\[
			\boldsymbol{r_{cs_b}}\times\boldsymbol{\ddot{r}_{cs_b}} =  0  \]
			\item The orbital radius is of a specified magnitude,
			\[ ||\boldsymbol{r_{cs_b}}|| =   d  \]
		\end{itemize}
	Based on these requirements, the following candidate Lyapunov function is selected
		\begin{equation}
		\begin{split}
		V(\boldsymbol{r_{cs_b}}, \boldsymbol{\dot{r}_{cs_b}}) = & \frac{1}{2}
			(|\boldsymbol{r_{cs_b}}\cdot\boldsymbol{\dot{r}_{cs_b}}|^2 + ||\boldsymbol{r_{cs_b}}\times\boldsymbol{\dot{r}_{cs_b}}
			 - \boldsymbol{L_d}||^2) \\ &
			+ \frac{a}{2}(||\boldsymbol{r_{cs_b}}|| - d)^2
		\end{split}			
			\end{equation}
		where $\boldsymbol{L_d}$ is the desired angular momentum and $a > 0 $ is a tuning parameter for 
		control design purpose.
		Now
\begin{equation*}
\begin{split}
\frac{dV}{dt} = & (\boldsymbol{r_{cs_b}}\cdot\boldsymbol{\dot{r}_{cs_b}})(||\boldsymbol{\dot{r}_{cs_b}}||^2 + \boldsymbol{r_{cs_b}}\cdot\boldsymbol{\ddot{r}_{cs_b}}) \\ & + (\boldsymbol{r_{cs_b}}\times\boldsymbol{\dot{r}_{cs_b}} - \boldsymbol{L_d})\cdot(\boldsymbol{r_{cs_b}}\times\boldsymbol{\ddot{r}_{cs_b}}) \\&
+a(||\boldsymbol{r_{cs_b}}|| - d)\frac{\boldsymbol{r_{cs_b}}\cdot\boldsymbol{\dot{r}_{cs_b}}}{||\boldsymbol{r_{cs_b}}||}
\end{split}
\end{equation*}
\noindent
Using the identity:
$(a\times b)\cdot(c\times d) = (a\cdot c)(b\cdot d) - (a\cdot d)(b\cdot c)$
\begin{equation*}
\begin{split}
(\boldsymbol{r_{cs_b}}\times\boldsymbol{\dot{r}_{cs_b}})\cdot(\boldsymbol{r_{cs_b}}\times\boldsymbol{\ddot{r}_{cs_b}}) = & ||\boldsymbol{r_{cs_b}}||^{2}(\boldsymbol{\dot{r}_{cs_b}}\cdot\boldsymbol{\ddot{r}_{cs_b}}) \\& - (\boldsymbol{r_{cs_b}}\cdot\boldsymbol{\ddot{r}_{cs_b}})(\boldsymbol{r_{cs_b}}\cdot\boldsymbol{\dot{r}_{cs_b}})
\end{split}
\end{equation*}
\noindent
and
\begin{equation}
\begin{split}
\frac{dV}{dt} = &(\boldsymbol{r_{cs_b}}\cdot\boldsymbol{\dot{r}_{cs_b}})||\boldsymbol{\dot{r}_{cs_b}}||^{2} + (\boldsymbol{r_{cs_b}}\cdot\boldsymbol{\ddot{r}_{cs_b}})(\boldsymbol{r_{cs_b}}\cdot\boldsymbol{\dot{r}_{cs_b}}) \\ &+ ||\boldsymbol{r_{cs_b}}||^{2}(\boldsymbol{\dot{r}_{cs_b}}\cdot\boldsymbol{\ddot{r}_{cs_b}}) - (\boldsymbol{r_{cs_b}}\cdot\boldsymbol{\ddot{r}_{cs_b}})(\boldsymbol{r_{cs_b}}\cdot\boldsymbol{\dot{r}_{cs_b}}) \\ &- \boldsymbol{L_d}\cdot(\boldsymbol{r_{cs_b}}\times\boldsymbol{\ddot{r}_{cs_b}})
+a(||\boldsymbol{r_{cs_b}}|| - d)\frac{\boldsymbol{r_{cs_b}}\cdot\boldsymbol{\dot{r}_{cs_b}}}{||\boldsymbol{r_{cs_b}}||} \\
\\
= &(\boldsymbol{r_{cs_b}}\cdot\dot{r}{cs_b})||\boldsymbol{\dot{r}_{cs_b}}||^{2} 
+ ||\boldsymbol{r_{cs_b}}||^{2}(\boldsymbol{\dot{r}_{cs_b}}\cdot\boldsymbol{\ddot{r}_{cs_b}}) \\ &- \boldsymbol{L_d}\cdot(\boldsymbol{r_{cs_b}}\times\boldsymbol{\ddot{r}_{cs_b}}) 
+a(||\boldsymbol{r_{cs_b}}|| - d)\frac{\boldsymbol{r_{cs_b}}\cdot\boldsymbol{\dot{r}_{cs_b}}}{||\boldsymbol{r_{cs_b}}||}
\end{split} \label{dVdt}
\end{equation}
\\
 Now we substitute (\ref{acc}) in the above equation and evaluate each term of the right hand side seperately to obtain

\begin{equation}
\begin{split}
\boldsymbol{\dot{r}_{cs_b}}\cdot\boldsymbol{\ddot{r}_{cs_b}} = &
\boldsymbol{\dot{r}_{cs_b}}\cdot \boldsymbol{r_{cs_b}}(\dot{\phi}^{2} - \frac{km_1}{||\boldsymbol{r_{1s_b}}||^3} - \frac{km_2}{||\boldsymbol{r_{2s_b}}||^3}) \\ & + \boldsymbol{\dot{r}_{cs_b}}\cdot(\dot{\phi}^2\boldsymbol{r_{c_b}} - \frac{km_1}{||\boldsymbol{r_{1s_b}}||^3}(\boldsymbol{r_{1c_b}}) \\ & - \frac{km_2}{||\boldsymbol{r_{2s_b}}||^3}(\boldsymbol{r_{2c_b}}))
+ \boldsymbol{\dot{r}_{cs_b}}\cdot\boldsymbol{\bar{u}_b}
 \end{split} \label{term_1}
\end{equation}

\begin{equation}
\begin{split}
\boldsymbol{r_{cs_b}}\times\boldsymbol{\ddot{r}_{cs_b}} 
= & \boldsymbol{r_{cs_b}}\times(\dot{\phi}^2\boldsymbol{r_{c_b}} - 2\hat{\boldsymbol{\omega}}\boldsymbol{\dot{r}_{cs_b}} - \frac{km_1}{||\boldsymbol{r_{1s_b}}||^3}(\boldsymbol{r_{1c_b}}) \\ & - \frac{km_2}{||\boldsymbol{r_{2s_b}}||^3}(\boldsymbol{r_{2c_b}}))  + \boldsymbol{r_{cs_b}}\times\boldsymbol{\bar{u}_b}
\end{split} \label{term_2}
\end{equation}
\noindent
\\
Evaluating (\ref{naturaldyn}) at L4 we get the following identity
\begin{equation}
\dot{\phi}^2\boldsymbol{r_{c_b}} - \frac{km_1}{||\boldsymbol{r_{1c_b}}||^3}(\boldsymbol{r_{1c_b}}) - \frac{km_2}{||\boldsymbol{r_{2c_b}}||^3}(\boldsymbol{r_{2c_b}}) = 0
\label{atL4}
\end{equation}
\noindent
	To simplify further analysis, we now make the following assumption:
	\begin{assumption}
	The radius of the orbit of the spacecraft is small compared to the distance of the libration point 
	from the two primaries such that $||\boldsymbol{r_{1s_b}}|| \approx ||\boldsymbol{r_{1c_b}}||(1 + \epsilon_1)$ and $||\boldsymbol{r_{2s_b}}|| \approx ||\boldsymbol{r_{2c_b}}||(1 + \epsilon_2)$ where $|\epsilon_i| << 1$. 
	\end{assumption}
	$\Box.$
	\\

	\noindent 
	With the above assumption and using a binomial identity, (\ref{term_1}) and (\ref{term_2}) reduce to

\begin{equation*}
\begin{split}
\boldsymbol{\dot{r}_{cs_b}}\cdot\boldsymbol{\ddot{r}_{cs_b}} =
& \boldsymbol{\dot{r}_{cs_b}}\cdot \boldsymbol{r_{cs_b}}(\dot{\phi}^{2} - \frac{km_1}{||\boldsymbol{r_{1s_b}}||^3} - \frac{km_2}{||\boldsymbol{r_{2s_b}}||^3}) \\
& + \boldsymbol{\dot{r}_{cs_b}}\cdot(\frac{km_13\epsilon_1}{||\boldsymbol{r_{1c_b}}||^3}\boldsymbol{r_{1c_b}} + \frac{km_23\epsilon_2}{||\boldsymbol{r_{2c_b}}||^3}\boldsymbol{r_{2c_b}}) \\
& + \boldsymbol{\dot{r}_{cs_b}}\cdot\boldsymbol{\bar{u}_b} 
 \end{split}
\end{equation*}

\begin{equation*}
\begin{split}
\boldsymbol{r_{cs_b}}\times\boldsymbol{\ddot{r}_{cs_b}}
= & \boldsymbol{r_{cs_b}}\times(-2\hat{\boldsymbol{\omega}}\boldsymbol{\dot{r}_{cs_b}})  	+ \boldsymbol{r_{cs_b}}\times\boldsymbol{\bar{u}_b} \\
& + \boldsymbol{r_{cs_b}}\times(\frac{km_13\epsilon_1}{||\boldsymbol{r_{1c_b}}||^3}\boldsymbol{r_{1c_b}} + \frac{km_23\epsilon_2}{||\boldsymbol{r_{2c_b}}||^3}\boldsymbol{r_{2c_b}})
\end{split}
\end{equation*}
Define
\begin{equation}
\boldsymbol{z} \deff  \frac{km_13\epsilon_1}{||\boldsymbol{r_{1c_b}}||^3}\boldsymbol{r_{1c_b}} + \frac{km_23\epsilon_2}{||\boldsymbol{r_{2c_b}}||^3}\boldsymbol{r_{2c_b}}
\end{equation}

\noindent
Thus (\ref{dVdt}) reduces to 
\begin{equation*}
\begin{split}
\frac{dV}{dt} = & ||\boldsymbol{r_{cs_b}}||^2(\boldsymbol{r_{cs_b}}\cdot\boldsymbol{\dot{r}_{cs_b}})(\dot{\phi}^2 + \frac{||\boldsymbol{\dot{r}_{cs_b}}||^2}{||\boldsymbol{r_{cs_b}}||^2} - \frac{km_1}{||\boldsymbol{r_{1s_b}}||^3} \\ &
- \frac{km_2}{||\boldsymbol{r_{2s_b}}||^3}) + \boldsymbol{z}\cdot(\boldsymbol{\dot{r}_{cs_b}}||\boldsymbol{r_{cs_b}}||^2 - \boldsymbol{L_d}\times \boldsymbol{r_{cs_b}})
 \\ & - \boldsymbol{L_d}\cdot(\boldsymbol{r_{cs_b}}\times(-2\hat{\boldsymbol{\omega}}\boldsymbol{\dot{r}_{cs_b}}))  + ||\boldsymbol{r_{cs_b}}||^{2}(\boldsymbol{\dot{r}_{cs_b}}\cdot\boldsymbol{\bar{u}_b}) \\ & - \boldsymbol{L_d}\cdot(\boldsymbol{r_{cs_b}}\times\boldsymbol{\bar{u}_b})
+a(||\boldsymbol{r_{cs_b}}|| - d)\frac{\boldsymbol{r_{cs_b}}\cdot\boldsymbol{\dot{r}_{cs_b}}}{||\boldsymbol{r_{cs_b}}||}
\end{split}
\end{equation*}
\noindent
Using the identity: $a\cdot(b\times c) = c\cdot({a\times b})$

\begin{equation}
\begin{split}
\frac{dV}{dt} = & ||\boldsymbol{r_{cs_b}}||^2(\boldsymbol{r_{cs_b}}\cdot\boldsymbol{\dot{r}_{cs_b}})(\dot{\phi}^2 + \frac{||\boldsymbol{\dot{r}_{cs_b}}||^2}{||\boldsymbol{r_{cs_b}}||^2} - \frac{km_1}{||\boldsymbol{r_{1s_b}}||^3} \\ &- \frac{km_2}{||\boldsymbol{r_{2s_b}}||^3}) -
(-2\hat{\boldsymbol{\omega}}\boldsymbol{\dot{r}_{cs_b}})\cdot(\boldsymbol{L_d}\times \boldsymbol{r_{cs_b}}) \\
&  + \boldsymbol{z}\cdot(\boldsymbol{\dot{r}_{cs_b}}||\boldsymbol{r_{cs_b}}||^2 - \boldsymbol{L_d}\times \boldsymbol{r_{cs_b}}) \\
& + \boldsymbol{\bar{u}_b}\cdot(\boldsymbol{\dot{r}_{cs_b}}||\boldsymbol{r_{cs_b}}||^2 - \boldsymbol{L_d}\times \boldsymbol{r_{cs_b}}) \\ &
+a(||\boldsymbol{r_{cs_b}}|| - d)\frac{\boldsymbol{r_{cs_b}}\cdot\boldsymbol{\dot{r}_{cs_b}}}{||\boldsymbol{r_{cs_b}}||}
\end{split} \label{dVdt2}
\end{equation}
\noindent
	Based on the need to render the time-derivative of the Lyapunov function to be
	 negative semi-definite, we select the controller as 
	\begin{equation}
	\begin{split}
		\boldsymbol{\bar{u}_b} = & -\beta(\boldsymbol{\dot{r}_{cs_b}}||\boldsymbol{r_{cs_b}}||^2 - \boldsymbol{L_d}\times \boldsymbol{r_{cs_b}}) +
			 \boldsymbol{\dot{r}_{cs_b}}\times\eta \\ &+ q\boldsymbol{r_{cs_b}} 
		- \frac{a(||\boldsymbol{r_{cs_b}}|| - d)}{||\boldsymbol{r_{cs_b}}||^3}\boldsymbol{r_{cs_b}} 
		- \boldsymbol{z}
	\end{split} \label{controller}
	\end{equation}
	where $\beta > 0$, $\eta \in \R^3$ and $q \in \R$ are to be chosen on further analysis. 
	Note that the controller is such that the first term renders the fourth term of (\ref{dVdt2}) as 
		negative. However the other terms do remain. Hence the additional two degrees of
		 freedom $\eta$ and $q$ appear in the proposed form.  

	Define
		\begin{equation}
		p \deff \dot{\phi}^2 + \frac{||\boldsymbol{\dot{r}_{cs_b}}||^2}{||\boldsymbol{r_{cs_b}}||^2} -
		 \frac{km_1}{||\boldsymbol{r_{1s_b}}||^3} - \frac{km_2}{||\boldsymbol{r_{2s_b}}||^3} 
		\end{equation}
	and rewrite 
\begin{equation*}
\begin{split}
\frac{dV}{dt} =& -\beta||(\boldsymbol{\dot{r}_{cs_b}}||\boldsymbol{r_{cs_b}}||^2 - \boldsymbol{L_d}\times \boldsymbol{r_{cs_b}})||^2 \\&
+ ||\boldsymbol{r_{cs_b}}||^2(\boldsymbol{r_{cs_b}}\cdot\boldsymbol{\dot{r}_{cs_b}})(p+q)\\ &
-(-2\hat{\boldsymbol{\omega}}\boldsymbol{\dot{r}_{cs_b}} + \boldsymbol{\dot{r}_{cs_b}}\times\eta)\cdot(\boldsymbol{L_d}\times \boldsymbol{r_{cs_b}}) \\ &+ (\boldsymbol{\dot{r}_{cs_b}}\times\eta)\cdot(\boldsymbol{\dot{r}_{cs_b}}||\boldsymbol{r_{cs_b}}||^2) 
- (q\boldsymbol{r_{cs_b}})\cdot(\boldsymbol{L_d}\times \boldsymbol{r_{cs_b}})
\end{split}
\end{equation*}

From the fact that 
\begin{equation*}
(\boldsymbol{\dot{r}_{cs_b}}\times\eta)\cdot(\boldsymbol{\dot{r}_{cs_b}}||\boldsymbol{r_{cs_b}}||^2) = 0
\end{equation*}
\begin{equation*}
(q\boldsymbol{r_{cs_b}})\cdot(\boldsymbol{L_d}\times \boldsymbol{r_{cs_b}}) = 0
\end{equation*}
\noindent
and choosing
\begin{equation}
q = -p \;\;\; \hbox{and}   \;\;\;  \eta = -2 \boldsymbol{\omega}
\end{equation}
that ensures 
\begin{equation*}
p + q = 0
\end{equation*}
\begin{equation*}
-2\hat{\boldsymbol{\omega}}\boldsymbol{\dot{r}_{cs_b}} + \boldsymbol{\dot{r}_{cs_b}}\times\eta = 0
\end{equation*}
we have 
\begin{equation}
\frac{dV}{dt} = -\beta||(\boldsymbol{\dot{r}_{cs_b}}||\boldsymbol{r_{cs_b}}||^2 - \boldsymbol{L_d}\times \boldsymbol{r_{cs_b}})||^2  \leq 0
\label{dV}
\end{equation}
\noindent
Therefore, the feedback system is stable in the sense of Lyapunov. Now we proceed to prove asymptotic stability of the desired orbit. 
\section{Application of LaSalle's invariance principle}
To prove asymptotic stability of the desired orbit, we employ LaSalle's invariance principle \cite{khalil2002nonlinear}. 	
	\begin{claim}
	Consider the positively invariant and compact set
		\begin{equation}
			\Omega =\{ (\boldsymbol{r_{cs_b}},\boldsymbol{\dot{r}_{cs_b}}) 
				\in (\R^3 \times \R^3) \hspace{1pt}|\hspace{1pt} V(\boldsymbol{r_{cs_b}},\dot{\boldsymbol{r_{cs_b}}}) \leq c \}
		\end{equation}
		The largest invariant set  in
		\[
		E = \{(\boldsymbol{r_{cs_b}},\boldsymbol{\dot{r}_{cs_b}}) \in \Omega | \dot{V} = 0\}
		\]
		is  $\{ (\boldsymbol{r_{cs_b}},\boldsymbol{\dot{r}_{cs_b}}): ||\boldsymbol{r_{cs_b}}|| = d, \boldsymbol{r_{cs_b}}\times\boldsymbol{\dot{r}_{cs_b}} = \boldsymbol{L_d} \}.$
	\end{claim}
	\proof 
From (\ref{dV}) it is clear that the set $\Omega$ is positively invariant. Also, from the construction of the set, it is evident that the set is closed and bounded. 
Since $\Omega \subset \mathbb{R}^{3\times 3}$, closed and bounded implies compactness. 
\\ Now consider the set $E = \{(\boldsymbol{r_{cs_b}},\boldsymbol{\dot{r}_{cs_b}}) \in \Omega | \dot{V} = 0\}$. This implies (from (\ref{dV}))
\begin{equation*}
-\beta||\boldsymbol{\dot{r}_{cs_b}}||\boldsymbol{r_{cs_b}}||^2 - \boldsymbol{L_d}\times \boldsymbol{r_{cs_b}}||^2 = 0
\end{equation*}
\noindent
which is
\begin{equation}
\boldsymbol{\dot{r}_{cs_b}}||\boldsymbol{r_{cs_b}}||^2 - \boldsymbol{L_d}\times \boldsymbol{r_{cs_b}} = 0 \label{LaSalles1}
\end{equation}
Taking the dot product of the above equation with $\boldsymbol{r_{cs_b}}$ yields
\begin{equation*}
\boldsymbol{r_{cs_b}}\cdot(\boldsymbol{\dot{r}_{cs_b}}||\boldsymbol{r_{cs_b}}||^2 - \boldsymbol{L_d}\times \boldsymbol{r_{cs_b}}) = 0
\end{equation*}
\begin{equation*}
(\boldsymbol{r_{cs_b}}\cdot\boldsymbol{\dot{r}_{cs_b}})||\boldsymbol{r_{cs_b}}||^2 - \boldsymbol{r_{cs_b}}\cdot(\boldsymbol{L_d}\times \boldsymbol{r_{cs_b}}) = 0
\end{equation*}
The second term vanishes, which gives
\begin{equation}
\boldsymbol{r_{cs_b}}\cdot\boldsymbol{\dot{r}_{cs_b}} = 0 \label{perpendicular}
\end{equation}
Now, taking the dot product of equation (\ref{LaSalles1}) with $\boldsymbol{\dot{r}_{cs_b}}$ yields
\begin{equation*}
\boldsymbol{\dot{r}_{cs_b}}\cdot(\boldsymbol{\dot{r}_{cs_b}}||\boldsymbol{r_{cs_b}}||^2 - \boldsymbol{L_d}\times \boldsymbol{r_{cs_b}}) = 0
\end{equation*}
\begin{equation*}
(\boldsymbol{\dot{r}_{cs_b}}\cdot\boldsymbol{\dot{r}_{cs_b}})||\boldsymbol{r_{cs_b}}||^2 - \boldsymbol{\dot{r}_{cs_b}}\cdot(\boldsymbol{L_d}\times \boldsymbol{r_{cs_b}}) = 0
\end{equation*}
\begin{equation}
||\boldsymbol{\dot{r}_{cs_b}}||^2||\boldsymbol{r_{cs_b}}||^2 - \boldsymbol{L_d}\cdot(\boldsymbol{r_{cs_b}}\times\boldsymbol{\dot{r}_{cs_b}}) = 0 \label{LaSalles2}
\end{equation}
Using equation (\ref{perpendicular}) we can write 
\begin{equation}
\boldsymbol{r_{cs_b}}\times\boldsymbol{\dot{r}_{cs_b}} = ||\boldsymbol{r_{cs_b}}||||\boldsymbol{\dot{r}_{cs_b}}||
\end{equation}
Equation (\ref{LaSalles2}) can be rewritten as 
\begin{equation*}
(\boldsymbol{r_{cs_b}}\times\boldsymbol{\dot{r}_{cs_b}})\cdot(\boldsymbol{r_{cs_b}}\times\boldsymbol{\dot{r}_{cs_b}}) - \boldsymbol{L_d}\cdot(\boldsymbol{r_{cs_b}}\times\boldsymbol{\dot{r}_{cs_b}}) = 0
\end{equation*}
\begin{equation}
\boldsymbol{L_d} = \boldsymbol{r_{cs_b}}\times\boldsymbol{\dot{r}_{cs_b}} \label{angMomntm}
\end{equation}
Equations (\ref{perpendicular}) and (\ref{angMomntm}) imply that the orbits lying in the limit set E are circular with constant angular momentum.
\\
\noindent
Further, the acceleration term (\ref{acc}) can be reduced as follows:
\begin{equation*}
\begin{split}
\boldsymbol{\ddot{r}_{cs_b}} = &\dot{\phi}^2\boldsymbol{r_{c_b}} - \frac{km_1}{||\boldsymbol{r_{1s_b}}||^3}\boldsymbol{r_{1c_b}} - \frac{km_2}{||\boldsymbol{r_{2s_b}}||^3}\boldsymbol{r_{2c_b}} - \frac{||\boldsymbol{\dot{r}_{cs_b}}||^2}{||\boldsymbol{r_{cs_b}}||^2}\boldsymbol{r_{cs_b}}
\\ &-\beta(\boldsymbol{\dot{r}_{cs_b}}||\boldsymbol{r_{cs_b}}||^2 - \boldsymbol{L_d}\times \boldsymbol{r_{cs_b}}) - \frac{a(||\boldsymbol{r_{cs_b}}|| - d)}{||\boldsymbol{r_{cs_b}}||^3}\boldsymbol{r_{cs_b}} \\ &
- (\frac{km_13\epsilon_1}{||\boldsymbol{r_{1c_b}}||^3}\boldsymbol{r_{1c_b}} + \frac{km_23\epsilon_2}{||\boldsymbol{r_{2c_b}}||^3}\boldsymbol{r_{2c_b}})
\end{split}
\end{equation*}
Using the assumption made earlier ($||\boldsymbol{r_{is_b}}|| \approx ||\boldsymbol{r_{ic_b}}||(1 + \epsilon_1)$ with $\epsilon_i << 1$) and equation (\ref{atL4}) we have 
\begin{equation*}
\begin{split}
\boldsymbol{\ddot{r}_{cs_b}} =& -\beta(\boldsymbol{\dot{r}_{cs_b}}||\boldsymbol{r_{cs_b}}||^2 - \boldsymbol{L_d}\times \boldsymbol{r_{cs_b}}) - \frac{||\boldsymbol{\dot{r}_{cs_b}}||^2}{||\boldsymbol{r_{cs_b}}||^2}\boldsymbol{r_{cs_b}} \\& - \frac{a(||\boldsymbol{r_{cs_b}}|| - d)}{||\boldsymbol{r_{cs_b}}||^3}\boldsymbol{r_{cs_b}}
\end{split}
\end{equation*}
When the system approaches steady state the first term in above equation renders zero. So we have following
\begin{equation}
\boldsymbol{\ddot{r}_{cs_b}} = - \frac{||\boldsymbol{\dot{r}_{cs_b}}||^2}{||\boldsymbol{r_{cs_b}}||^2}\boldsymbol{r_{cs_b}} - \frac{a(||\boldsymbol{r_{cs_b}}|| - d)}{||\boldsymbol{r_{cs_b}}||^3}\boldsymbol{r_{cs_b}} \label{reducedacc}
\end{equation}
Since the steady state orbit is a circular orbit with constant speed and angular momentum, the acceleration above should be equal to the centrepetal acceleration 
\begin{equation}
\boldsymbol{\ddot{r}_{cs_b}} = - \frac{||\boldsymbol{\dot{r}_{cs_b}}||^2}{||\boldsymbol{r_{cs_b}}||^2}\boldsymbol{r_{cs_b}} \label{centrepetal}
\end{equation}
Equations (\ref{reducedacc}) and (\ref{centrepetal}) imply second term in (\ref{reducedacc}) should be zero, which results in
\begin{equation}
||\boldsymbol{r_{cs_b}}|| = d \label{radius}
\end{equation}
Therefore the steady state orbit is a circular orbit (\ref{perpendicular}) with constant angular momentum (\ref{angMomntm}) and desired radius (\ref{radius}). 
\par
\remark The control law (\ref{controller}) is observed to have a feedback linearization structure. It can be written as
\begin{equation*}
\begin{split}
\boldsymbol{\bar{u}_b} = & -\beta\boldsymbol{e_1} - a\boldsymbol{e_2} - f(\boldsymbol{r_{cs_b}}, \boldsymbol{\dot{r}_{cs_b}}) - \frac{||\boldsymbol{\dot{r}_{cs_b}}||^2}{||\boldsymbol{r_{cs_b}}||^2}\boldsymbol{r_{cs_b}} 
\end{split}
\end{equation*}
where $e_1 = \boldsymbol{\dot{r}_{cs_b}}||\boldsymbol{r_{cs_b}}||^2 - \boldsymbol{L_d}\times \boldsymbol{r_{cs_b}}$ and $e_2 = \frac{(||\boldsymbol{r_{cs_b}}|| - d)}{||\boldsymbol{r_{cs_b}}||^3}\boldsymbol{r_{cs_b}}$ represent error in the current state and desired state, $f(\boldsymbol{r_{cs_b}}, \boldsymbol{\dot{r}_{cs_b}})$ is the natural dynamics of the spacecraft given by the right hand side of (\ref{naturaldyn}) and the fourth term is the dersired dynamics.

\begin{figure*}[!htb]
\begin{minipage}[t][][b]{0.3\textwidth}
  \includegraphics[width=\linewidth]{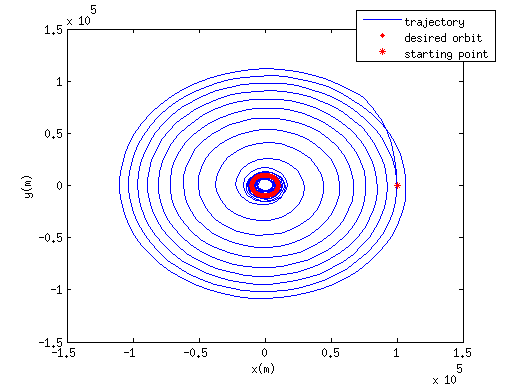}
  \caption{Case1: Spaceraft tajectory}
  \label{traj2D}
\end{minipage}\hspace{0.1cm}
\begin{minipage}[t][][b]{0.3\textwidth}
  \includegraphics[width=\linewidth]{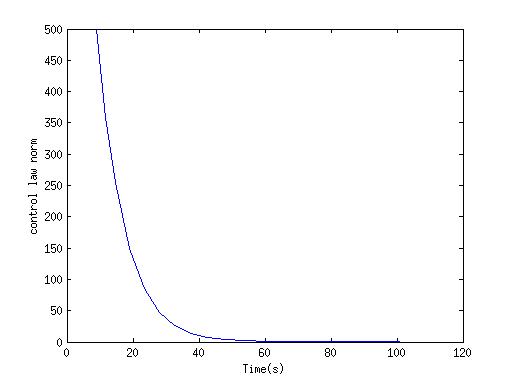}
  \caption{Case1: Norm of the control - initial few seconds}
  \label{u12D}
\end{minipage}\hspace{0.1cm}
\begin{minipage}[t][][b]{0.3\textwidth}
  \includegraphics[width=\linewidth]{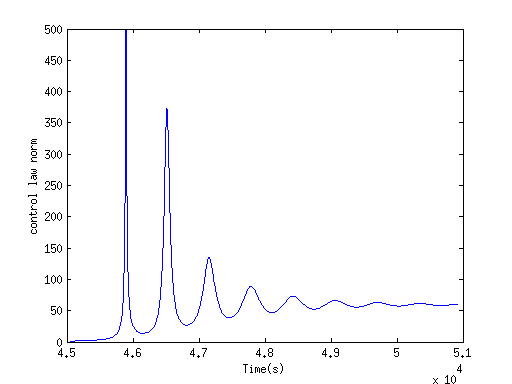}
  \caption{Case1: Norm of the control - just before and after orbit capture}
  \label{u22D}
\end{minipage}

\end{figure*}

\begin{figure*}[!htb]
\begin{minipage}{0.3\textwidth}
  \includegraphics[width=\linewidth]{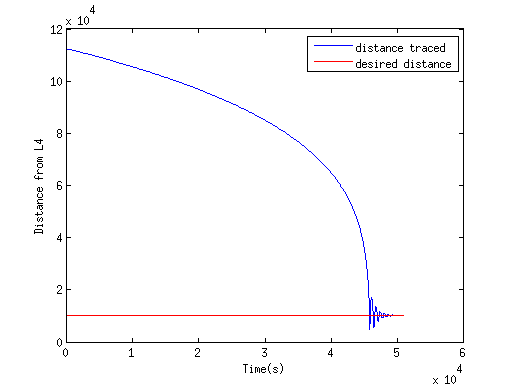}
  \caption{Case1: Distance of spacecraft from L4}
  \label{r2D}
\end{minipage}
\begin{minipage}{0.3\textwidth}
  \includegraphics[width=\linewidth]{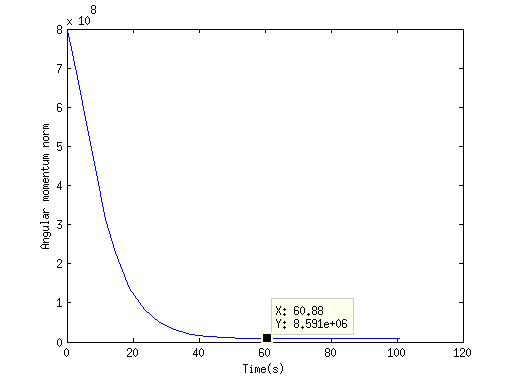}
  \caption{Case1: Norm of angular momentum - initial few seconds}
  \label{L12D}
\end{minipage}
\begin{minipage}{0.3\textwidth}
  \includegraphics[width=\linewidth]{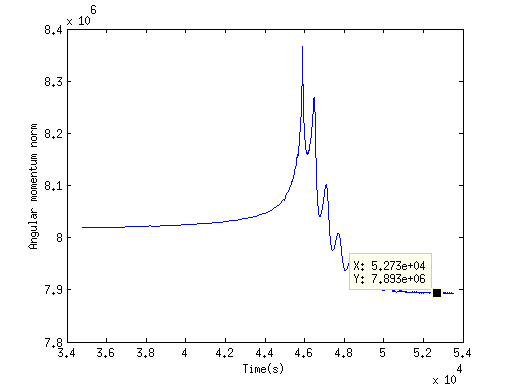}
  \caption{Case1: Norm of angular momentum - just before and after orbit capture}
  \label{L22D}
\end{minipage}
\end{figure*}

\section{Numerical simulations}

Simulations were performed for a spacecraft initially located near the Earth-Moon L4 point using a numerical integration scheme. For computing the control law it is required to know the distance of the spacecraft, $||\boldsymbol{r_{1s_b}}||$ and $||\boldsymbol{r_{2s_b}}||$, from both the primaries. 
To keep the control magnitude under bounds, a saturation was introduced and the control 
law implemented was
\begin{equation}
\boldsymbol{u} = \begin{cases}
\boldsymbol{\bar{u}_b} &\text{if $||\boldsymbol{\bar{u}_b}|| \leq u_{max}$}\\
u_{max}\frac{\boldsymbol{\bar{u}_b}}{||\boldsymbol{\bar{u}_b}||} &\text{if $||\boldsymbol{\bar{u}_b}|| > u_{max}$}
\end{cases}
\end{equation}
where $\boldsymbol{\bar{u}_b}$ is given by (\ref{controller}). Two representative 
cases are presented here. The first one corresponds to the initial position in the plane of rotation of the two primaries. The second one corresponds to an arbitrary initial position in the 3D space. Parameters and conditions corresponding to both the cases are given below. 

\begin{figure*}[!htb]
\begin{minipage}[t][][b]{0.3\textwidth}
  \includegraphics[width=\linewidth]{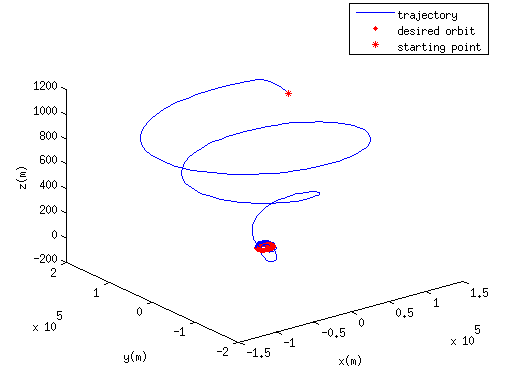}
  \caption{Case2: Spaceraft tajectory}
  \label{traj3D}
\end{minipage}\hspace{0.1cm}
\begin{minipage}[t][][b]{0.3\textwidth}
  \includegraphics[width=\linewidth]{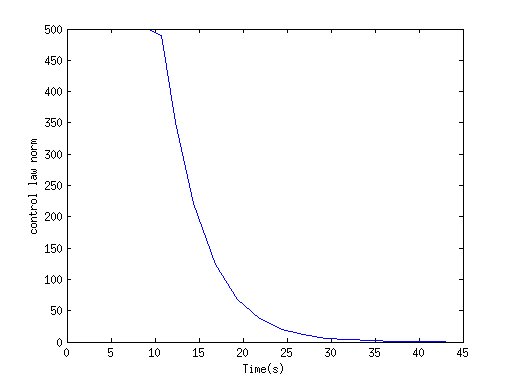}
  \caption{Case2: Norm of control - initial few seconds }
  \label{u13D}
\end{minipage}\hspace{0.1cm}
\begin{minipage}[t][][b]{0.3\textwidth}
  \includegraphics[width=\linewidth]{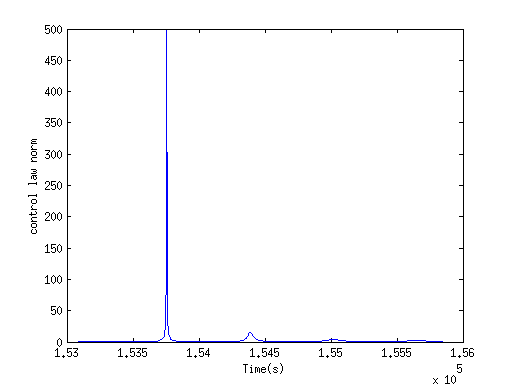}
  \caption{Case2: Norm of the control - just before and after orbit capture}
  \label{u23D}
\end{minipage}

\end{figure*}

\begin{figure*}[!htb]
\begin{minipage}{0.3\textwidth}
  \includegraphics[width=\linewidth]{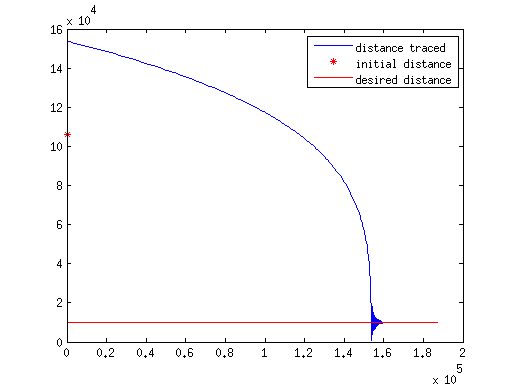}
  \caption{Case2: Distance of spacecraft from L4}
  \label{r3D}
\end{minipage}
\begin{minipage}{0.3\textwidth}
  \includegraphics[width=\linewidth]{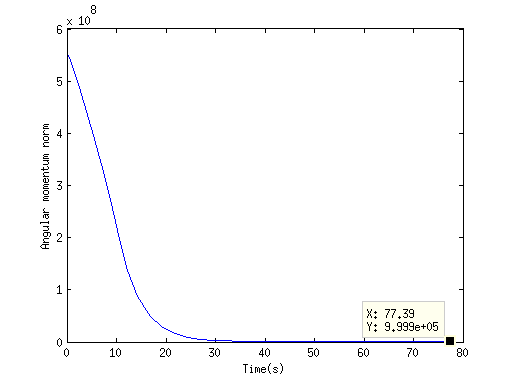}
  \caption{Case2: Norm of angular momentum - initial few seconds}
  \label{L13D}
\end{minipage}
\begin{minipage}{0.3\textwidth}
  \includegraphics[width=\linewidth]{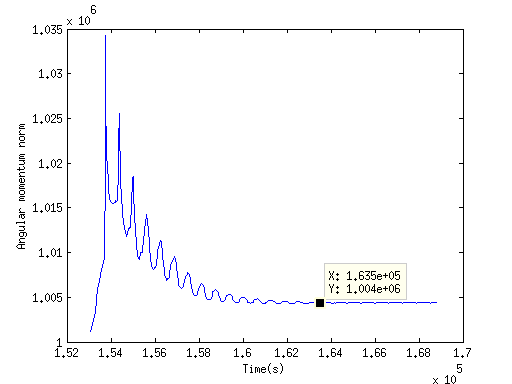}
  \caption{Case2: Norm of angular momentum - just before and after orbit capture}
  \label{L23D}
\end{minipage}
\end{figure*}

\begin{table}[H]
\centering
\begin{tabular}{|l|l|l|l|}
\hline
\textbf{Parameter} & \textbf{Case 1 (2D)} & \textbf{Case 2 (3D)} & \textbf{Units} \\ \hline 
$\boldsymbol{r_{cs_b}}$ at t = 0 & [100000 0 0] & [75000 75000 1000] & m \\ \hline
 $\boldsymbol{\dot{r}_{cs_b}}$ at t = 0 & [0 8000 0] & [100 7500 10] & m$s^{-1}$ \\ \hline
 $\boldsymbol{L_d}$ & [0 0 80000000] & [0 0 1000000] & $m^2$ $s^{-1}$\\ \hline
 $d$ & 10000 & 10000 & m\\ \hline
 $u_{max}$ & \multicolumn{2}{|c|}{500} & $m$ $s^{-2}$  \\ \hline
 $\beta$ & \multicolumn{2}{|c|}{1e-11} & $m^{-2}$ $s^{-1}$ \\ \hline
$a$ & \multicolumn{2}{|c|}{$d^2/10000$ } &  $m_2$ $s^{-2}$ \\ \hline
  $\boldsymbol{\omega}$ & \multicolumn{2}{|c|}{[0 0 2.66e-06]} & $rad$ $s^{-1}$ \\ \hline
$k$ & \multicolumn{2}{|c|}{6.673e-11} & $N$ $m^2$ $kg^{-2}$ \\ \hline
$m_1$ & \multicolumn{2}{|c|}{5.972e24}& $kg$ \\
Earth's mass & \multicolumn{2}{|c|}{} & \\ \hline
$m_2$ & \multicolumn{2}{|c|}{7.34767e22} & $kg$ \\ 
Moon's mass & \multicolumn{2}{|c|}{} & \\ \hline
\end{tabular}
\caption{Parameters and conditions for the two cases}
\end{table}

Figure(\ref{traj2D}) shows the trajectory of the spacecraft in the first case. It can be seen that it  takes about 13.5 hours to converge to the desired orbit. In the initial few seconds the control law shoots to a maximum and then becomes almost zero (fig(\ref{u12D})). It remains 
close to zero till the spacecraft reaches the desired orbit. Near the desired orbit there is a spike in the control input. The spacecraft distance from the L4 point 
oscillates about the desired distance and then finally settles down to a circular orbit with radius within 0.5\% of the desired radius. This shows that the initial control input provides the spacecraft a momentum such that it goes into a spiral orbit. The sudden spike in the control near desired orbit is similar to orbit insertion burns performed by rocket thrusters in spacecraft to enter an orbit around a planetary body.
\par
The results for the first case of 2D motion are similar to the second case of 3D motion. The trajectory followed by the spacecraft depends on the initial error as it will contribute to the momentum imparted intially. The larger the error, larger is the momentum imparted and hence can lead to larger oscillations about the desired orbit before it stabilizes. Further, in the stable state a continous control input is required to maintain the spacecraft in the orbit. This is the centrepetal force and depends on the desired angular momentum.  

\section{Conclusion and Future work}
In this work, using Lyapunov theory we have shown that asymptotically stable orbits around the triangular libration point can be achieved with low thrust solutions. The initial thrust is required for a short time to transfer from an arbitrary intial position to a transfer orbit. Subsequent stationkeeping requires negligible thrust. This concept can be extended to transfer from one circular orbit to another. 
\par 
This analysis did not include the perturbation effects by solar gravity and other planets. Further work could include these factors and modify the control law accordingly. Also the assumption about $\epsilon_i$ requires that we determine a region of stability around the L4 point. This can be done by observing the unstabe manifolds of the collinear libration points. This can then lay the ground work for spacecraft formations. 

\addtolength{\textheight}{-12cm}   



\bibliography{references}{}

\end{document}